\documentstyle[12pt,dina4]{article}
 
\newcommand{\subsetc}{\subset \! \! \! \! \raisebox{-3.5pt}{$\to $} \, }
\newtheorem{THEOREM}{Theorem}[section]
\newcommand{\thm}[1]{
\begin{THEOREM}
#1
\end{THEOREM}
}
\newtheorem{PROPOSITION}[THEOREM]{Proposition}
\newcommand{\prop}[1]{
\begin{PROPOSITION}
#1
\end{PROPOSITION}
}
\newtheorem{rem}{Remark}[section]

\newtheorem{LEMMA}[THEOREM]{Lemma}
\newcommand{\la}[1]{
\begin{LEMMA}
#1
\end{LEMMA}
} 
\newtheorem{COROLLARY}[THEOREM]{Corollary}

\newtheorem{DEFINITION}{Definition}[section]
 
\newcommand{\pf}{
\begin{flushleft}
{\bf Proof}
\end{flushleft}}
\newcommand{\nm}[1]{\parallel #1 \parallel}

\newcommand{\nmwg}[1]{\nm{#1}_{L^p_{v(\gamma )}}}

\begin{document}
 \begin{center}
 {\LARGE{}On a restriction problem of de Leeuw type\\ for Laguerre multipliers
 }\\[.4cm]
 {\sc{}George  Gasper\footnote{Department of Mathematics, Northwestern 
 University, Evanston, IL 60208, USA. The work of this author was 
 supported in part by the National Science   Foundation under grant 
 DMS-9103177.}
 and Walter  Trebels\footnote{Fachbereich Mathematik, TH Darmstadt, 
Schlo\ss{}gartenstr.7, D--64289
 Darmstadt, Germany.} }\\[.4cm]
 {\it Dedicated to K\'aroly Tandori on the occasion of his 70th birthday}\\
[.4cm]
 {(Aug. 10, 1994 version)}\\[.4cm]
 \end{center}

 \bigskip
 {\bf Abstract.} 
In 1965 K. de Leeuw \cite{deleeuw} proved among other things
in the Fourier transform setting: 
{\it If a continuous function $m(\xi _1, \ldots ,\xi _n)$ on ${\bf R}^n$ 
generates a bounded transformation on $L^p({\bf R}^n),\; 1\le p \le \infty ,$ 
then its trace  
$\tilde{m}(\xi _1, \ldots ,\xi _m)=m(\xi _1, \ldots ,\xi _m,0,\ldots ,0), \; 
m<n,$ generates a bounded transformation on $L^p({\bf R}^m)$. }
In this paper, the analogous problem is discussed in the setting of 
Laguerre expansions of different orders. 


\bigskip
{\bf AMS(MOS) subject classifications.} 33C45, 42A45, 42C10 

\bigskip 
\section{Introduction}

The purpose of this paper is to discuss the question: suppose 
$\{ m_k\} _{k\in {\bf N}_0}$ generates a bounded transformation with respect to a Laguerre 
function expansion of order $\alpha $ on some $L^p$--space, does it also 
generate a corresponding bounded map with respect to a Laguerre function 
expansion  of order $\beta \, $? To become more precise let us 
first introduce some notation. 
Consider the  Lebesgue spaces
$$L^p_{w(\gamma )} = \{ f: \; \nmwg{f} = ( \int _0^\infty
|f(x) e^{-x/2}|^p x^{\gamma }\, dx )^{1/p} < \infty \} \; ,\quad 1 \le p <
\infty  ,$$
$$L^\infty _{w(\gamma )} =\{ f: \; \| f\|_{L^\infty _{w(\gamma )}} = {\rm ess 
\; sup} _{x>0}|f(x)e^{-x/2}| < \infty \}, \; \quad p=\infty ,$$
where $\gamma >-1$. Let $L_n^\alpha (x), \; \alpha >-1,\; n\in {\bf N}_0$,
denote the classical Laguerre polynomials (see Szeg\"o 
\cite[p. 100]{szego})  
and set 
$$ R_n^\alpha (x)=L_n^\alpha (x)/L_n^\alpha (0),\quad \; \; \;  L_n^\alpha (0)=
A_n^\alpha =
\left( \matrix{ n+\alpha \cr n} \right) =\frac{\Gamma (n+ \alpha +1)}{\Gamma
(n+1) \Gamma (\alpha +1)  }.$$
Associate to $f$ its formal Laguerre series
$$ f(x) \sim (\Gamma (\alpha +1))^{-1}
\sum _{k=0}^\infty {\hat f}_\alpha (k) L_k^\alpha (x), $$
where the Fourier Laguerre coefficients of $f$ are defined by
\begin{equation}\label{coeff}
{\hat f}_\alpha (n) = \int _0^\infty f(x)R_n^\alpha (x)x^\alpha e^{-x}\, dx
\end{equation}
(if the integrals exist). 
A sequence $m=\{ m_k\} _{k\in {\bf N}_0} $ is called a (bounded) multiplier 
from $L^p_{w(\gamma )}$ into $L^q_{w(\delta )}$, notation 
$ m \in M^{p,q}_{\alpha ;\, \gamma ,\, \delta }$, if
$$ \| \sum _{k=0}^\infty m_k{\hat f}_\alpha (k) L_k^\alpha \| _{L^q
_{w(\delta )}} \le C 
\| \sum _{k=0}^\infty {\hat f}_\alpha (k) L_k^\alpha \| _{L^p_{w(\gamma )}} $$
for all polynomials $f$; the smallest constant $C$ for which this holds is 
called the multiplier norm $ \| m\| _{M^{p,q}_{\alpha ;\, \gamma ,\, \delta }}. $ 
For the sake of simplicity we write $M^{p,q}_{\alpha ;\, \gamma }:= M^{p,q} 
_{\alpha ;\, \gamma ,\, \gamma }$ if $\gamma =\delta $ and, 
if additionally $p=q$, 
$M^p_{\alpha ; \, \gamma }:=M^{p,p}_{\alpha ;\, \gamma }$.

\bigskip \noindent
We are {\bf mainly} interested in the question: when is 
$M^{p,q}_{\alpha ;\, \alpha }$ continuously embedded in 
$M^{p,q}_{\beta ;\, \beta }$:
$$M^{p,q}_{\alpha ;\, \alpha } \subsetc 
M^{p,q}_{\beta ;\, \beta }\, ,  \quad 1\le p\le q \le \infty \; ?$$
The Plancherel theory immediately yields 
$$ l^\infty =M^2_{\alpha ; \alpha } =M^2_{\beta ; \beta },\quad 
\alpha ,\beta >-1.$$ 

\smallskip \noindent
A combination of sufficient multiplier conditions with necessary ones indicates 
which results are to be expected. To this end, define the  
 fractional difference operator $\Delta ^\delta $ of order $\delta$ by
$$ \Delta ^\delta m_k = \sum _{j=0}^\infty A_j^{-\delta -1}m_{k+j}$$
(whenever the series converges), the classes $wbv_{q,\delta }\, ,\; 1\le q \le 
\infty ,\;  \delta >0, $ of weak bounded variation (see \cite{wbv}) of  
bounded sequences which have finite norm $\| m \| _{q,\delta }\, ,$ where 
\begin{eqnarray*}
 \| m\| _{q,\delta } := \sup _k |m_k| + \sup _{N\in {\bf N}_0} 
\left( \sum _{k=N}^{2N} |(k+1)^\delta \Delta ^\delta m_k|^q \frac{1}{k+1} 
\right) ^{1/q}, &  & q <\infty , \\
 \| m\| _{\infty ,\delta } := \sup _k |m_k| + \sup _{N\in {\bf N}_0} 
|(k+1)^\delta \Delta ^\delta m_k| \quad \quad \quad \quad \quad \quad \quad ,
&  & q =\infty .
\end{eqnarray*}
Observing the duality (see \cite{sttr})
\begin{equation}\label{dual}
 M^p_{\alpha ;\, \gamma }=M^{p'}_{\alpha ;\, \alpha p' -\gamma p'/p}\, ,
\quad -1 < \gamma < p(\alpha +1) -1, \quad 1<p<\infty ,
\end{equation}
where $1/p+1/p'=1$, 
we may restrict ourselves to the case $1<p<2$. The Corollary 1.2 b) in 
\cite{sttr} gives the embedding 
\begin{equation}\label{nec}
 M^p_{\alpha ;\, \alpha } \subsetc wbv_{p',s},\quad s=(2\alpha +2/3)(1/p-1/2),\quad 
\alpha >-1/3,
\end{equation}
when $(2\alpha +2)(1/p-1/2) >1/2$. Theorem 5 in \cite{wbv} gives the first 
embedding in 
$$ wbv_{p',s} \subsetc wbv_{2,s} \subsetc M^p_{\beta ;\, \beta},$$
whereas the last one follows from Corollaries 1.2 and 4.5 
in \cite{sttr} provided $s> {\rm max} \{ (2\beta +2)(1/p-1/2),\; 1\} ,\; \beta 
>-1.$ Hence, choosing $\gamma =\alpha $ in (\ref{dual}), we obtain
\prop{ 
Let $1<p<\infty $ and $\alpha $ be such that  $(2\alpha +2/3)|1/p-1/2|>1$. 
Then
$$M^p_{\alpha ;\, \alpha } \subsetc M^p_{\beta ;\, \beta }\, ,
\quad -1<\beta < \alpha -2/3.$$
}

\medskip \noindent 
In the same way we can derive a result for $M^{p,q}$--multipliers.
The necessary condition in \cite[Cor. 1.3]{laguerre} can easily be extended 
in the sense of \cite[Cor. 2.5 b)]{laguerre} to
$$ \sup _k |(k+1)^\sigma m_k| +\sup _n (\sum _{k=n}^{2n} |(k+1)^{\sigma +s} 
\Delta ^sm_k|^{q'}/k)^{1/q'} \le C \| m\| _{M^{p,q}_{\alpha ;\, \alpha }}\, ,$$
where $\alpha >-1/3,\; 1/q=1/p-\sigma /(\alpha +1),\; 1<p<q<2 ,
\; (\alpha +1)(1/q-1/2)>1/4,$ and $s=(2\alpha +2/3)(1/q-1/2)>0 $. Using this 
and the sufficient condition for $M^{p,q}_{\beta ;\, \beta }$--multipliers 
given in \cite[Cor. 1.2]{gst}, which is proved only for $\beta \ge 0$, 
we obtain 
$$ M^{p,q}_{\alpha ;\, \alpha } \subsetc M^{p,q}_{\beta ;\, \beta }\, ,
\quad 0\le \beta < \alpha -2/3, \; (2\alpha +2/3)(1/q-1/2)>1,\;  
1<p<q<2.$$  

\medskip \noindent
In this context let us mention that the same technique yields for $1<p,q<2$
\begin{equation}\label{p,q}
M^p_{\alpha ;\, \alpha } \subsetc M^q_{\beta ;\, \beta }\, ,\quad (2\alpha +2/3) 
( 1/p-1/2) > {\rm max}\{ ( 2\beta +2)( 1/q-1/2),\; 1\}  .
\end{equation}
This embedding is in so far interesting as it allows to go from $p,\, 1<p<2,$ 
to $q\neq p,\, 1<q<2,$ connected with a loss in 
the size of $\beta $ if $q <p$ or a gain in $\beta $ if $1<p<q<2\, $; 
e.g. 
$$M^{4/3}_{10;\, 10} \subsetc M^q_{5;\, 5}\, ,\quad 1.08 \le q \le 2,\quad 
{\rm or} \quad M^{8/7}_{2;\, 2} \subsetc M^q_{4;\, 4}\, ,\quad 3/2\le q \le 2.$$

\smallskip \noindent
Improvements of (\ref{p,q}) can be expected by better necessary conditions 
and/or better sufficient conditions; but this technique {\bf cannot} give 
something like 
$$M^p_{\alpha;\, \alpha } \subsetc 
M^q_{\beta ;\, \beta },\quad (\alpha +1)(1/p-1/2) >(\beta +1)(1/q-1/2),\; 
1<p,\, q<2,$$
which is suggested by (\ref{p,q}) when choosing ``large'' $\alpha $ with $p$ 
near $2$ since then the number
$4(1/p-1/2)/3$, which describes the smoothness gap between the necessary 
conditions and the sufficient conditions in \cite[Cor. 1.2]{sttr}, 
is ``negligible''.

\bigskip \noindent
Concerning the general problem ``When does $M^{p,q}_{\alpha ;\, \gamma _1,\, 
\delta _1} \subsetc M^{p,q}_{\beta;\, \gamma _2,\, \delta _2} $ 
hold?'',\\
we mention results in Stempak and Trebels \cite[Cor. 4.3]{sttr}: For 
$1<p<\infty $ there holds 
$$ M^p_{\beta ; \beta p/2 + \delta } =
M^p _{0; \delta } \quad {\rm if} \; \left\{ \begin{array}{r@{\; ,\quad}l}
-1-\beta p/2 < \delta <p-1+\beta p/2 & -1<\beta <0 \, ,\\
-1<\delta < p-1 & 0 \le \beta ,
\end{array} \right.  $$
which for $\delta =0$ contains Kanjin's \cite{kan} result 
and for $\delta =p/4-1/2$ Thangavelu's \cite{thanga}.
In particular, there holds for $-1 < \beta <\alpha ,\; 1<p<\infty $, 
\begin{equation}\label{transmul}
 M^p_{\beta ; \, \beta } = M^p_{\beta ;\, \beta p/2 +\beta p(1/p-1/2)} =
M^p_{\alpha ;\,  \alpha p/2 + \beta p(1/p-1/2)}\, ,\quad 
(2\beta +2)|1/p-1/2| <1.
\end{equation}
These results are based on Kanjin's \cite{kan} transplantation theorem and 
its weighted version in \cite{sttr}. The latter gives further insight 
into our problem in so far as it implies that the restriction $\beta < \alpha 
-2/3$ in Proposition 1.1 is not sharp. \\
To this end we first note that the following extension of Corollary 4.4 in 
\cite{sttr} holds
$$ wbv_{2,s} \subsetc M^p_{\alpha ;\, \alpha p/2 +\eta (p/2-1)},\quad 0\le \eta 
\le 1,\quad 1<p\le 2,\quad s>1/p.$$
(For the proof observe that for $\alpha =0$ the parameter $\gamma =\eta 
(p/2-1),\; 0\le \eta \le 1,$ is admissible in \cite[Theorem 1.1]{sttr} and then 
follow the argumentation of \cite[Cor. 4.4]{sttr}.) This combined with
(\ref{nec}) yields for $s=(2\alpha +2/3)(1/p-1/2)>1/p$
$$ M^p_{\alpha ;\, \alpha } \subsetc wbv_{2,s} \subsetc 
M^p_{\alpha ; \, \alpha p/2 +p/2 -1 }\, \, , \quad 1<p\le 2, \quad \alpha 
>(p+1)/(6-3p). $$
Thus, by interpolation with change of measure,
$$M^p_{\alpha ;\, \alpha } \subsetc M^p_{\alpha ;\, \alpha p/2+\delta }\, ,
\quad p/2-1 \le \delta \le \alpha -\alpha p/2,\quad \alpha 
>(p+1)/(6-3p). $$
Since (\ref{transmul}) gives 
$$M^p_{\alpha ;\, \alpha p/2 +\beta 
p(1/p-1/2)} = M^p_{\beta ;\, \beta } $$
we arrive at
\prop{
Let $1<p\le 2$ and $\alpha >(p+1)/(6-3p) $. Then
$$ M^p_{\alpha ;\, \alpha } \subsetc M^p_{\beta ;\, \beta }\, ,
\quad (2\beta +2)(1/p-1/2)<1,\quad -1<\beta < \alpha  .$$
}

\medskip \noindent
The first restriction on $\beta $ is equivalent to $\beta <(2p-2)/(2-p)$. This 
combined with the restriction on $\alpha $ gives $\alpha -\beta 
>(7-5p)/(6-3p)$,  
the latter being decreasing in $p$ and taking the value $2/3$ at $p=1$. Hence 
Proposition 1.2 is an improvement of the previous one for all $1<p<2$ provided 
$(p+1)/(6-3p)<\alpha \le (2p-2)/(2-p).$ For big $\alpha $'s, Proposition 1.1 is 
certainly better.
If in the transplantation theorem in \cite{sttr} higher exponents could be 
allowed in the power weight -- which is possible in the Jacobi expansion case 
as shown by Muckenhoupt \cite{mu} -- the technique 
just used would give the embedding when $-1<\beta < \alpha ,\; 1<p<2,$ and
$\alpha >(p+1)/(6-3p).$\\
Summarizing, it is reasonable to 
\begin{flushleft}
{\large\bf conjecture} \hspace{2cm} 
$M^{p,q}_{\alpha;\, \alpha } \subsetc 
M^{p,q}_{\beta ;\, \beta },\quad -1<\beta <\alpha ,\quad 1\le p \le q \le 
\infty .$
\end{flushleft}
Apart from the above fragmentary results, so far we can only 
prove the conjecture in the extreme case when $q=\infty $ and $\beta 
\ge 0$; the latter restriction arises from the  fact that we have to make use 
of the twisted Laguerre convolution (see \cite{indag}) which is proved till now 
only for Laguerre polynomials $L_n^\alpha (x)$ with $\alpha \ge 0$. 
Our main result is
\thm{
If $1\le p \le \infty $, then 
$$M^{p,\infty }_{\alpha ;\, \alpha } \subsetc
M^{p,\infty }_{\beta ;\, \beta },\quad 0\le \beta <\alpha .$$
}

\bigskip \noindent
{\bf Remarks.} 1) One could speculate that an interpolation argument applied to
$$M^2_{\alpha ;\, \alpha }=M^2_{\beta ;\, \beta }\, ,\quad 
M^\infty _{\alpha ;\, \alpha }=M^1_{\alpha ;\, \alpha } \subsetc
M^1_{\beta ;\, \beta } =M^\infty _{\beta ;\, \beta }\, ,\quad \beta < \alpha ,$$
could give the open case $M^p_{\alpha ;\, \alpha } \subsetc M^p_{\beta ;\, 
\beta }\, ,\; 1<p<2$. In this respect we mention a result of Zafran 
\cite[p. 1412]{za} for the Fourier transform pointed out to us by A. Seeger:
\begin{quote}
{\it Denote by $M^p({\bf R}) $ the set of bounded Fourier multipliers on 
$L^p({\bf R})$ and by $M^\wedge ({\bf R})$ the set of Fourier transforms of 
bounded measures on ${\bf R}$. Then $M^p({\bf R}),\; 1<p<2, $ 
is {\bf not} an interpolation space with respect to the 
pair $(M^\wedge ({\bf R}) ,L^\infty ({\bf R}) )$. }
\end{quote}
Thus de Leeuw's result mentioned at the beginning cannot be proved by 
interpolation.

\medskip \noindent
2) It is perhaps amazing to note that the $wbv$--classes do not play only an 
auxiliary role in dealing with the above formulated general problem. 
In the framework of one-dimensional Fourier transforms/series this was shown by
Muckenhoupt, Wheeden, and Wo-Sang Young \cite{mwy}. That this 
phenomenon also occurs in the framework of Laguerre expansions can be seen 
from the following two theorems.
\thm{
If $\alpha >-1,\; \alpha \neq 0 $, then 
$$ wbv _{2,1} \subsetc M^2_{\alpha ;\, \alpha +1} \, .$$
In the case $-1 <\alpha <0$ the multiplier operator is defined only on the 
subspace $\{ f\in L^2_{w(\alpha +1)}: {\hat f}_\alpha (0) =0 \} $. 
}
\thm{
If $\alpha > -1$, then 
$$ M^2_{\alpha ;\, \alpha +1} \subsetc wbv _{2,1}\, .$$
}
A combination of these two results leads to 
\begin{equation}\label{char}
 M^2_{\alpha ;\, \alpha +1} = M^2_{\beta ;\, \beta +1} =wbv_{2,1}\, ,\quad 
\alpha ,\beta > -1,\quad \alpha ,\beta \neq 0,
\end{equation} 
and a combination with \cite[(19)]{sttr} gives
$$ M^2_{\alpha ;\, \alpha +1} \subsetc M^p_{\alpha ;\, \alpha }\, ,\quad 
\alpha \ge 0,\quad (2\alpha +2)/(\alpha +1) <p\le 2.$$

\bigskip
\section{Proof of Theorem 1.3} 
Theorem 1.3 is an immediate consequence of the combination of the following two 
theorems. 

\medskip \noindent 
\thm{
Let $f\in L^p_{w(\alpha )}$ with $\alpha >-1 $ when $ 1\le p < \infty $ and 
$\alpha \ge 0$ when $p=\infty $. Then there 
exists a function $g\in L^p_{w(\beta )},\; -1< \beta < \alpha ,$ with
$$ g(x) \sim (\Gamma (\beta +1))^{-1}
\sum _{k=0}^\infty {\hat f}_\alpha (k) L_k^\beta (x),\quad 
\| g\| _{L^p_{w(\beta )}} \le C \| f\| _{L^p_{w(\alpha )}}\, .$$
}
\pf{ 
First let $1\le p < \infty$ and, without loss of generality, let $f$ be 
a polynomial (these are dense in $L^p_{w(\alpha )}$). We recall
the projection formula (3.31) in Askey and Fitch \cite{af}
$$ e^{-x}L_n^\beta (x) =\frac{1}{\Gamma (\alpha -\beta )} \int _x^\infty 
(y-x)^{\alpha -\beta -1}e^{-y}L_n^\alpha (y)\, dy, \quad -1<\beta <\alpha . $$
Then the following computations are justified.
\begin{eqnarray*}
\| g\| _{L^p_{w(\beta )}} & = & C \left( \int _0^\infty |\sum _{k=0}^\infty 
{\hat f}_\alpha (k) L_k^\beta (x)e^{-x/2}|^p x^\beta dx \right) ^{1/p} \\
& = & C\left( \int _0^\infty \left| \int _x ^\infty (y-x)^{\alpha -\beta -1} 
e^{-y} \sum _{k=0}^\infty {\hat f}_\alpha (k) L_k^\alpha (y)\, dy \right| ^p 
x^\beta e^{xp/2}dx \right) ^{1/p} \\
& \le & C \int _1^\infty (t-1)^{\alpha -\beta -1}\left( \int _0^\infty | \sum 
_k {\hat f}_\alpha (k) L_k^\alpha (xt) x^{\alpha -\beta +\beta /p } 
e^{-x(t-1/2)}|^p dx\right) ^{1/p}dt 
\end{eqnarray*}
after a substitution and application of the integral Minkowski inequality. 
Additional substitutions lead to
\begin{eqnarray*}
\| g\| _{L^p_{w(\beta )}} & \le & C \int _0^\infty s^{\alpha -\beta -1}
(s+1)^{\beta /p' -\alpha -1/p} \times \\
 &  & \quad \quad \left( \int _0^\infty | \sum 
_k {\hat f}_\alpha (k) L_k^\alpha (y) e^{-y/2} y^{(\alpha -\beta)/p' } 
e^{-ys/2(s+1)}|^p y^\alpha dy \right) ^{1/p}ds \\ 
& \le & C \int _0^\infty s^{(\alpha -\beta )/p -1}
(s+1)^{-(\alpha +1)/p }
\left( \int _0^\infty | \sum 
_k {\hat f}_\alpha (k) L_k^\alpha (y) e^{-y/2} |^p y^\alpha dy\right) ^{1/p}ds, 
\end{eqnarray*}
where we used the inequality $y^{(\alpha -\beta)/p' } e^{-ys/2(s+1)} \le C 
((s+1)/s)^{(\alpha -\beta )/p'}$. Since $-1<\beta <\alpha $ it is easily seen 
that the outer integration only gives a bounded contribution.

\smallskip \noindent
If $f\in L^\infty _{w(\alpha )}$ then $|(k+1)^{-1/2}{\hat f}_\alpha (k)| 
\le C \| f\| _{L^\infty _{w(\alpha )}}$ by \cite[Lemma 1]{analy} and, 
therefore, the Abel-Poisson means of an arbitrary $f\in L^\infty _{w(\alpha )}$ can be 
represented by 
$$ P_rf(x) =(\Gamma (\alpha +1))^{-1}
\sum _k r^k {\hat f}_\alpha (k) L_k^\alpha (x),\quad 0\le r<1,\quad x\ge 0,$$
and, by the convolution theorem in G\"orlich and Markett \cite[p. 169]{indag},
$$\| P_rf\| _{L^\infty _{w(\alpha )}} \le C \| f\| _{L^\infty _{w(\alpha )}}, 
\quad 0\le r<1,\quad \alpha \ge 0.$$ 
A slight modification of the argument in the case $1\le p < \infty $ shows that 
$$ \| g_r\| _{L^\infty _{w(\beta )}}:= \| (\Gamma (\beta +1))^{-1}
\sum _k r^k {\hat f}_\alpha (k) L_k^\beta \| _{L^\infty _{w(\beta )}} \le 
C \| P_rf\| _{L^\infty _{w(\alpha )}} \le C\| f\| _{L^\infty _{w(\alpha )}} .$$
By the weak${}^*$ compactness there exists a function $g\in L^\infty 
_{w(\beta )}$ with ${\hat g}_\beta (k) = {\hat f}_\alpha (k)  $ and $\| g\| 
_{L^\infty _{w(\beta )}} \le \lim \inf _{k\to \infty } \| g_{r_k}\| 
_{L^\infty _{w(\beta )}}$ for a suitable sequence $r_k\to 1^-$; hence also the 
assertion in the case $p=\infty $.  
}
\thm{
For $\alpha \ge 0$ there holds

\smallskip \noindent
i) \hspace{1.8cm}  $M^{1,p}_{\alpha ;\, \alpha } =M^{p',\infty}_{\alpha ;\, 
\alpha } = L^p_{w(\alpha )}, \quad 1< p\le \infty $,

\smallskip \noindent
ii) \hspace{1.7cm} $M^{1,1}_{\alpha ;\, \alpha } =M^{\infty ,\infty}_{\alpha ;\, 
\alpha } =\{ m=\{ m_k\} _{k\in {\bf N}_0}
: \| P_r(m)\| _{L^1_{w(\alpha )}} =O(1),\; r\to 1^- \} $, \\
where $P_r(m)(x) = (\Gamma (\alpha +1))^{-1} \sum_k r^km_kL_k^\alpha (x)$. 
}
\pf{
The first equalities in $i)$ and $ii)$ are the standard duality statements. Let us 
briefly indicate the second equalities (which are also more or less standard).

\smallskip \noindent
If $m=\{ m_k\} _{k\in {\bf N}_0}$ are the Fourier Laguerre coefficients 
of an $L^p_{w(\alpha )}$--function, $1<p\le \infty $, 
or in the case $p=1$ of a bounded measure 
with respect to the weight $e^{-x/2}x^\alpha $, then Young's inequality in 
G\"orlich and Markett \cite{indag} (or a slight extension of it to measures in 
the case $p=1$) shows that $m\in M^{p',\infty}_{\alpha ;\alpha }$.

\smallskip \noindent
Conversely, associate formally to a sequence $m=\{ m_k\}$ an operator $T_m$ by
\begin{equation}\label{op}
 T_mf(x) \sim (\Gamma (\alpha +1))^{-1}
\sum _{k=0}^\infty m_k {\hat f}_\alpha (k) L_k^\alpha (x).
\end{equation}
Then, in essentially the notation of G\"orlich and Markett \cite{indag},
$$ T_m(P_rf)(x)=P_r(m)*f(x)=  \int _0^\infty T_x^\alpha 
(P_r(m)(y))f(y)e^{-y}y^\alpha dy,$$
where $T_x^\alpha $ is the Laguerre translation operator.
If $\| f\| _{L^{p'}_{w(\alpha )}} =1$ then 
$$\| T_m(P_rf)\|_{L^\infty _{w(\alpha )}} \le \| m\| _{M^{p',\infty }_{\alpha 
;\, \alpha }} \| P_rf\| _{L^{p'}_{w(\alpha )}}
\le C \| m\| _{M^{p',\infty }_{\alpha ;\, \alpha }} \; ,$$
and hence, by the converse of H\"older's inequality,
$$ \sup _{\| f\| _{L^{p'}_{w(\alpha )}} =1} \left| \int _0^\infty 
T_x^\alpha (P_r(m)(y))e^{-y/2}y^{\alpha /p}f(y)e^{-y/2}y^{\alpha /p'}dy \right| 
\quad \quad \quad $$
$$ \quad \quad \quad \quad \quad =\| T_x^\alpha (P_r(m))\| _{L^p_{w(\alpha )}}
\le C \| m\| _{M^{p',\infty }_{\alpha ;\, \alpha }} $$
for $x\ge 0,\; 0\le r<1.$ In particular, for $x=0$ we obtain
$$\| P_r(m)\| _{L^p_{w(\alpha )}}
\le C \| m\| _{M^{p',\infty }_{\alpha ;\, \alpha }}\; ,\; \; 0\le r<1.$$
Now weak${}^*$ compactness gives the desired converse embedding.
}

\section{Proof of Theorems 1.4 and 1.5} 
The proof relies heavily on the Parseval formula
\begin{equation}\label{par}
\frac{1}{\Gamma (\alpha +1)} \sum _{k=0}^\infty A_k^\alpha |{\hat f}_\alpha 
(k)|^2 =\int _0^\infty |f(x)e^{-x/2}|^2 x^\alpha dx
\end{equation}
and its extension 
\begin{equation}\label{pardif}
\sum _{k=0}^\infty A_k^{\alpha +\lambda } |\Delta ^\lambda 
{\hat f}_\alpha (k)|^2 \approx 
\int _0^\infty |f(x)e^{-x/2}|^2 x^{\alpha +\lambda } dx,\quad \lambda \ge 0,
\end{equation}
which is a consequence of the formula 
\begin{equation}\label{lcoeff}
\Delta ^\lambda {\hat f}_\alpha (k)=
C_{\alpha ,\lambda } {\hat f}_{\alpha +\lambda }(k) 
\end{equation}
(see e.g. the proof of 
Lemma 2.1 in \cite{laguerre}).
For the proof of Theorem 1.4 we further need the following discrete 
analog of the $p=2$ case of a weighted Hardy inequality in Muckenhoupt 
\cite{muhardy} whose 
proof can at once be read off from \cite{muhardy} by replacing the integrals 
there by sums and using the fact that
$$a \le 2(a+b)^{1/2} [ (a+b)^{1/2} -b^{1/2} ]  $$
when $a,\, b \ge 0$; also see the extensions in \cite[Sec. 4]{anhe}.
\la{
Let  $\{ u_k\} _{k\in {\bf N}_0}\, ,\{ v_k\} _{k\in {\bf N}_0} $ be
non-negative sequences (if $v_k=0$ we set $v_k^{-1}=0$). Then 
\noindent
\begin{enumerate}
\item[a)]$ \quad \quad \quad \quad \displaystyle{ 
 \sum _{k=0}^\infty \Big| \sum _{j=0}^k a_j \Big| ^2u_k
\le C \sup _N \Big( \sum _{k=N}^\infty u_k \sum _{k=0}^N 
v_k^{-1} \Big) \sum _{j=0}^\infty | a_j | ^2v_j. }$
\item[b)]$ \quad \quad \quad \quad \displaystyle{ 
 \sum _{k=0}^\infty \Big| \sum _{j=k}^\infty a_j \Big| ^2u_k
\le C \sup _N \Big( \sum _{k=0}^N u_k \sum _{k=N}^\infty 
v_k^{-1} \Big)  \sum _{j=0}^\infty | a_j | ^2v_j. }$
\end{enumerate}
}

\smallskip

\bigskip \noindent
{\bf Proof of Theorem 1.4.} 
Using (\ref{pardif}) and the operator $T_m$ defined in (\ref{op}), we obtain 
$$ \int _0^\infty |T_mf(x)e^{-x/2}|^2 x^{\alpha +1}dx \approx 
\sum _{k=0}^\infty A_k^{\alpha +1} |\Delta (m_k{\hat f}_\alpha (k))|^2 .$$
Since 
\begin{equation}\label{firstdif}
\Delta (m_k{\hat f}_\alpha (k))=m_k\Delta {\hat f}_\alpha (k)
+{\hat f}_\alpha (k+1)\Delta m_k
\end{equation}
we first observe that 
$$ \sum _{k=0}^\infty A_k^{\alpha +1} |m_k|^2 |\Delta {\hat f}_\alpha (k)|^2 
\le \| m\| ^2_\infty \sum _{k=0}^\infty A_k^{\alpha +1} | \Delta {\hat f}
_\alpha (k)|^2 \le C \| m\| ^2_\infty \| f\| _{L^2_{w(\alpha +1)}} ^2 \, .$$
To dominate the term containing $\Delta m_k $ we deduce from (\ref{par}) that 
for $\alpha \ge 0 $ the Fourier Laguerre coefficients tend to zero as $
k\to \infty $. Hence
$$ \sum _{k=0}^\infty A_k^{\alpha +1} |{\hat f}_\alpha (k+1)\Delta m_k |^2
=\sum _{k=0}^\infty A_k^{\alpha +1} |\Delta m_k |^2\, 
\Big| \sum _{j=k+1}^\infty \Delta {\hat f}_\alpha (j)\Big| ^2=:I.$$
In order to apply Lemma 3.1 b), we choose $u_k=A_k^{\alpha +1} 
|\Delta m_k |^2$ and $v_k=A_k^{\alpha +1} $, and observe that when 
$M\in {\bf N},\, 2^{M-1} \le N< 2^M,$ we have that 
\begin{eqnarray*}
\Big( \sum _{k=0}^N u_k \sum _{k=N}^\infty 
v_k^{-1} \Big)  & 
 \le & C (N+1)^{-\alpha }\sum _{j=0}^M \sum _{k=2^j-1} 
^{2^{j+1}-2}(k+1)|\Delta m_k |^2 \frac{A_k^{\alpha +1}}{k+1}\\
& \le & C (N+1)^{-\alpha }\sum _{j=0}^M (2^j)^\alpha \| m\| ^2_{2,1}
\le C \| m\| ^2_{2,1}
\end{eqnarray*}
uniformly in $N$ if $\alpha >0$. Then Lemma 3.1 b) gives
$$ I \le C \| m\| _{2,1}^2 \sum _{j=0}^\infty A_j^{\alpha +1} 
|\Delta {\hat f}_\alpha (j)|^2 \le C\| m\| _{2,1}^2 \| f\| _{L^2_{w(\alpha 
+1)}}^2 $$
by (\ref{pardif}). Thus there remains to consider the case $-1 <\alpha <0$.
For the same choice of $u_k$ and $v_k$ one easily obtains 
$$ \Big( \sum _{k=N}^\infty  u_k \sum _{k=0}^N
v_k^{-1} \Big) \le C \| m\| _{2,1}^2 \, .$$
Now assume that ${\hat f}(0)=0$. Then we have
$$\sum _{k=0}^\infty A_k^{\alpha +1} |{\hat f}_\alpha (k+1)\Delta m_k |^2
=\sum _{k=0}^\infty A_k^{\alpha +1}|\Delta m_k|^2 \, \Big| \sum 
_{j=0}^k \Delta {\hat f}_\alpha (j)\Big| ^2 \le C\| m\| _{2,1}^2 \| f\| ^2
_{L^2_{w(\alpha +1)}} \, ,$$
where the last estimate follows by Lemma 4.1 a); thus Theorem 1.4 is 
established.

\smallskip \noindent

\bigskip \noindent
The {\bf proof of Theorem 1.5} is essentially contained in \cite{laguerre}. 
As in \cite{laguerre}, consider a monotone decreasing $C^\infty $-function 
$\phi(x)$ with 
$$ \phi (x) =
             \cases{$1$ & if $0 \le x \le 2$ \cr
                     $0$ & if $x\ge 4$
                   },
\quad \phi _i(x) = \phi (x/2^i).$$  
Then the $\phi _i(k)$ are the Fourier Laguerre coefficients of an 
$L^2_{w(\alpha +1)}$-function $\Phi ^{(i)}$ with norm $\| \Phi ^{(i)} \| 
_{L^2_{w(\alpha +1)}} \le C \, (2^i)^{\alpha /2}$ and
$$ \sum _{k=2^i}^{2^{i+1}}  A_k^{\alpha +1} | \Delta m_k|^2 
= \sum _{k=2^i}^{2^{i+1}} A_k^{\alpha +1} | \Delta (m_k\phi _i(k))|^2 
\le \sum _{k=0}^{2^{i+2}} A_k^{\alpha +1} | \Delta (m_k\phi _i(k))|^2 $$
$$\le C \| T_m\Phi ^{(i)} \| _{L^2_{w(\alpha +1)}}^2 
\le C \| m\| _{M^2_{\alpha ;\, \alpha +1}}^2 \| \Phi ^{(i)} \| 
_{L^2_{w(\alpha +1)}}^2  \le C 2^{i\alpha }
\| m\| _{M^2_{\alpha ;\, \alpha +1}}^2 .$$
This immediately leads to
$$ \| m\| _\infty  + \left( \sum _{2^i}^{2^{i+1}}  |(k+1)\Delta m_k|^2 
\frac{1}{k+1} \right) ^{1/2}
\le C \| m\| _{M^2_{\alpha ;\, \alpha +1}} ,$$
uniformly in $i$, since by \cite[(10)]{laguerre} there holds 
$\| m\| _\infty \le C\| m\| _{M^2_{\alpha ;\, \alpha +1}} $; thus Theorem 1.5 
is established.

\medskip \noindent
{\bf Remark.} 3) (Added on Aug. 10, 1994) The characterization (\ref{char}) 
can easily be extended to 
\begin{equation}\label{charex}
M^2_{\alpha ,\alpha +l} = wbv_{2,l}, \quad \alpha > -1, 
\quad \alpha \neq 0, \ldots , l-1,\; l\in {\bf N}.
\end{equation}
In the case $\alpha < l-1$ the multiplier operator is defined only on the 
subspace $\{ f\in L^2_{w(\alpha +l)}: {\hat f}_\alpha (k)=0,\; 0 \le k < 
(l-1-\alpha )/2 \} .$

\bigskip \noindent
The necessity part carries over immediately (see also \cite{laguerre}). The 
sufficiency part will be proved by induction. Thus suppose that (\ref{charex}) 
is true for $l=1, \ldots ,n $ and $\alpha $'s as indicated. Then, as in the 
case $n=1$, by (\ref{pardif})
$$ \int _0^\infty |T_mf(x)e^{-x/2}|^2 x^{\alpha +n+1}dx \approx 
\sum _{k=0}^\infty A_k^{\alpha +n+1} |\Delta ^n \Delta 
(m_k{\hat f}_\alpha (k))|^2 $$
$$\le C \sum _{k=0}^\infty A_k^{\alpha +n+1} |\Delta ^n (m_k \Delta 
{\hat f}_\alpha (k))|^2 + C \sum _{k=0}^\infty A_k^{\alpha +n+1} |\Delta ^n 
({\hat f}_\alpha (k+1)\Delta m_k )|^2 =: I+II $$
By the assumption and (\ref{lcoeff})
$$I \le C \| m\| ^2_{wbv_{2,n}}
\sum _{k=0}^\infty A_k^{\alpha +n+1} |\Delta ^n {\hat f}_{\alpha +1}(k))|^2 
 \le C \| m\| ^2_{wbv_{2,n+1}} \int _0^\infty |f(x)e^{-x/2}|^2 x^{\alpha 
+n+1} dx $$
on account of the embedding properties of the $wbv$--spaces \cite{wbv}. 
Analogously $II$ can be estimated by
$$II \le C \| \{ (k+1) \Delta m_k \} \| ^2_{wbv_{2,n}} \sum _{k=0}^\infty 
A_k^{\alpha +n+1} \Big| \Delta ^n \Big( \frac{{\hat f}_\alpha 
(k+1)}{k+1} \Big) \Big| ^2 .$$
By the Leibniz formula for differences there holds 
\begin{eqnarray*}
 \Delta ^n \Big( \frac{{\hat f}_\alpha (k+1)}{k+1} \Big) & \le & 
C \sum _{j=0}^n |\Delta ^j{\hat f}_\alpha (k+1) | \, |\Delta ^{n-j} 
\frac{1}{j+k+1} |   \\
{} & \le  & C \sum _{j=0}^n (j+k+1)^{j-n-1}|\Delta ^j{\hat f}_\alpha (k+1) | .
\end{eqnarray*}
Hence we have to dominate for $j=0,\, \ldots ,\, n$
$$ II_j:= \sum _{k=0}^\infty A_k^{\alpha -n-1+2j} | \Delta ^j{\hat f}_\alpha 
(k+1) |^2.$$
If $\alpha >n$ then $c_j:=-\alpha -2j+n+1 <1$ for all $j=0, \ldots ,n$, 
$\Delta ^j{\hat f}_\alpha (k+1) = \sum _{i=k+1} ^\infty \Delta ^{j+1}
{\hat f}_\alpha (i) $, 
and we can apply \cite[Theorem 346]{hlp} repeatedly to obtain 
$$II_j \le C \sum _{k=0}^\infty A_k^{\alpha -n-1+2j} |(k+1) \Delta ^{j+1}
{\hat f}_\alpha (k+1) |^2 \approx \sum _{k=0}^\infty A_k^{\alpha -n+2j+1} 
| \Delta ^{j+1}{\hat f}_\alpha (k+1) |^2$$
$$\le \ldots \le C \sum _{k=0}^\infty A_k^{\alpha +n+1} 
| \Delta ^{n+1}{\hat f}_\alpha (k+1) |^2 \le C 
\int _0^\infty |f(x)e^{-x/2}|^2 x^{\alpha +n+1} dx \, .$$
Since $\| \{ (k+1) \Delta m_k \} \| _{wbv_{2,n}} \le C \| m \| _{wbv_{2,n+1}} 
$, this gives the assertion for the weight $x^{n+1}$ in the case $\alpha > n.$

\medskip \noindent
If $\alpha <n,\; \alpha \neq 0, \ldots , n,$ then some $c_j>1$. For the 
application of \cite[Theorem 346]{hlp} one needs $c_j\neq 1$; this is 
guaranteed by the hypothesis $\alpha \neq 0,\dots , n$ (in the case of an 
additional weight $x^{n+1}$). For the $j$ for which $c_j>1$ 
we have to use the representation 
$$ \Delta ^j {\hat f}_\alpha (k+1) =-\sum _{i=0}^k \Delta ^{j+1} 
{\hat f}_\alpha (i),\quad {\rm if} \; \Delta ^j {\hat f}_\alpha (0)=0,$$
i.e., the first $(j+1)$ Fourier-Laguerre coefficients have to vanish to ensure 
this representation. But $0\le j \le j_0$, where $j_0$ is choosen in such a way 
that $c_{j_0}>1$ and $c_{j_0+1}<1$, hence $j_0 = [(n-\alpha )/2] $ (with 
respect to the additional weight $x^{n+1}$); here we used the standard notation 
for $[a],\; a\in {\bf R},$ to be the greatest integer $\le a$.
Hence the condition that the first $[(n-\alpha)/2]+1$ Fourier-Laguerre 
coefficients have to vanish is needed if the additional weight is $x^{n+1}$. 
A repeated application of \cite[Theorem 346]{hlp} with appropriate $c>1$ or 
$c<1$ now gives the assertion.

\end{document}